\documentstyle[11pt]{article}
\begin{document}


\title{On Explicit Formula for Restricted Partition Function}
\author{Boris Y. Rubinstein\\
Department of Engineering Sciences and Applied Mathematics,\\
Northwestern University, 2145 Sheridan Road \\
Evanston IL 60208-3125, U.S.A.}
\date{\today}

\maketitle

\begin{abstract}

A new recursive procedure of the calculation of a restricted partition
function is suggested. An explicit combinatorial formula for
the restricted partition function is found based on this procedure.


\end{abstract}


\section{Introduction}
The problem of partitions of positive integers has long history started from
the work of Euler \cite{Euler} who "laid a foundation of the theory of
partitions" \cite{GAndrews}, introducing the idea of generating functions.
Many great mathematicians, like Cayley, Sylvester, MacMahon, Ramanujan, and
others contributed to the development of the theory, using Euler idea.

Cayley \cite{Cayley} found explicit formulas for number $p_k(n)$ of partitions
of positive integer $n$ into at most $k$ parts with small $k$.
He also suggested a method of decomposition of the corresponding
generating function $$
  G_k(t) = \prod_{i=1}^k \frac{1}{1-t^i} =
  \sum_{n=0}^{\infty} p_k(n) \;t^n,
$$
and gave the combinatorial formula for such decompositions (unfortunately, this
formula itself requires knowledge of all partitions of $k$).

Sylvester was the next mathematician who provided a new insight and made a
remarkable progress in this field. He introduced \cite{Sylv0} the so-called
Ferrers graphs for presentation of partitions. He also found
\cite{Sylv1,Sylv2} the procedure
enabling to determine a {\it restricted} partition functions, and
described symmetry features of such functions. The restricted
partition function $p(n,{\bf d}^m) \equiv p(n,\{d_1,d_2,\ldots,d_m\})$ is a
number of partitions of $n$ into positive integers  $\{d_1,d_2,\ldots,d_m\}$,
each not greater than $n$. It is very simple to show that the generating
function in this case takes the form
\begin{equation}
G({\bf d}^m,t)=\prod_{i=1}^m\frac{1}{1-t^{d_{i}}}
 =\sum_{n=0}^{\infty} p(n,{\bf d}^m)\;t^n
\;.
\label{genfunc}
\end{equation}
Sylvester showed that the restricted partition
function may be presented as a sum of "waves",
each wave closely related to
prime roots of unit of degree $n$, where $n$ are prime divisors of elements
of the set ${\bf d}^m$. This fact was known to Herschel \cite{Herschel} who
introduced a notion of {\it circulator} and Cayley who used its elegant version
called {\it prime circulator} (see \cite{Gupta} for more information). Namely,
Sylvester showed that each wave $W_i$, where $i$ runs over distinct factors in
$d_1,d_2,...,d_m$, is a coefficient of ${t}^{-1}$ in the series
expansion in ascending powers of $t$ of
\begin{equation}
e^{s w_k}\;\;\prod_{r=1}^m \frac{1}
{1-e^{d_r u_k}}\;\;,\;\;w_k=2\pi i\;
\frac{p_k}{q}+t\;\;,\;\;u_k=2\pi i\; \frac{p_k}{q}-t\;,
\label{syl2}
\end{equation}
and $p_1,p_2,...,p_{max\;k}$ are integers (unity included) smaller than
$i$ and prime to it.
It should be noted here that the above result is
only a recipe for calculation of the partition function and doesn't provide an
explicit formula.

Sylvester found \cite{Sylv2} that the shifted partition function
$$
 q(n,{\bf d}^m) \equiv p(n-\frac{1}{2}\sum_{i=1}^m d_i,{\bf d}^m)
$$
has following parity properties:
$$
q(n,{\bf d}^{2m}) = -q(-n,{\bf d}^{2m}), \ \ \ \ \
q(n,{\bf d}^{2m+1}) = q(-n,{\bf d}^{2m+1}),
$$
and established that these functions have zeros at
all integer values of $n$ from $0$ to $m/2-1$ for even $m$ and
at all semiinteger values from $1/2$ to $m/2-1$ for odd $m$. He suggested
to use knowledge of partition function zeroes for its construction using the
method of indeterminate coefficients.

Recently, a different presentation
of the partition function, which may be called {\it polynomial expansion}, was
introduced in \cite{FelRub}, where a new recursive
procedure for calculation of the restricted partition function is found.
It permits to reconstruct in an unified way nearly all terms of expansion,
except one, which also demands usage of the
method of indeterminate coefficients.


In this article I present an approach enabling to overcome this
inconsistence, and to determine all terms of expansion in the framework of a
single recursive procedure. This method also produces the explicit formula of
the restricted  partition function.


\section{Polynomial part of the partition function}

J.J. Sylvester showed that the restricted partition function may be written as
a sum of "waves", he found the recipe for calculation of each such wave. We
consider in this section purely polynomial part of the partition function which
corresponds to the wave $W_1$. It may be found as a coefficient of $t^{-1}$ in an expansion of
the generator
\begin{equation}
G_1(s,t) = \frac{e^{st}}{\prod_{i=1}^m (1-e^{-d_i t})}.
\label{generator1}
\end{equation}
Sylvester found that $W_1(s)$ depends on Bernoulli numbers and sums of powers
of elements $d_i$.
We find an explicit form of polynomial part $W_1(s)$ through the Bernoulli
polynomials of higher order. We start from the generating function for the Bernoulli
polynomials of higher order $B_n^{(m)}(s | d_1,d_2,\ldots,d_m)$ \cite{bat53}:
\begin{equation}
\left(\prod_{i=1}^m d_i \right) t^m
\frac{e^{st}}{\prod_{i=1}^m (e^{d_i t}-1)} =
\sum_{n=0}^{\infty} B_n^{(m)}(s | {\bf d}^m) \frac{t^n}{n!},
\label{B_n^m_genfunc}
\end{equation}
where we use a shortcut notation
$$
{\bf d}^m = \{d_1,d_2, \ldots,d_m\}.
$$
One immediately obtains a presentation of $W_1$ generator in the form
\begin{equation}
G_1(s,t) =
\frac{1}{\pi({\bf d}^m) t^{m}}
\sum_{n=0}^{\infty} B_n^{(m)}(s | -{\bf d}^m) \frac{t^n}{n!},
\label{G_1}
\end{equation}
where
$$
\pi({\bf d}^m) =  \prod_{i=1}^m d_i.
$$
The coefficient of $1/t$ in the above expression is given by the term with
$n=m-1$
\begin{equation}
W_1(s,{\bf d}^m) =
\frac{1}{(m-1)!\pi({\bf d}^m)}
B_{m-1}^{(m)}(s | -{\bf d}^m).
\label{W_1}
\end{equation}
It is useful to consider a {\it shifted} restricted partition function defined
as follows
\begin{equation}
V(s,{\bf d}^m) = W(s-\xi ({\bf d}^m),{\bf d}^m)\;,
\ \ \xi ({\bf d}^m) = \frac{1}{2} \sum_{i=1}^m d_i\;,
\label{transform1}
\end{equation}
and its polynomial part $V_1(s,{\bf d}^m)$ which is cast in
\begin{equation}
V_1(s,{\bf d}^m) =
\frac{1}{(m-1)!\pi({\bf d}^m)}
B_{m-1}^{(m)}(s - \xi ({\bf d}^m) | -{\bf d}^m).
\label{V_1}
\end{equation}
Now we may use formula for Bernoulli polynomials of higher order found by
N\"orlund \cite{Norlund}
$$
B_n^{(m)}(s|-{\bf d}^m) = B_n^{(m)}(s+\sum_{i=1}^m d_i|{\bf d}^m)
$$
to arrive at
\begin{equation}
V_1(s,{\bf d}^m) =
\frac{1}{(m-1)!\pi({\bf d}^m)}
B_{m-1}^{(m)}(s + \xi ({\bf d}^m) | {\bf d}^m).
\label{V_1a}
\end{equation}
Then we apply another formula
$$
B_n^{(m)}(s|{\bf d}^m) = \sum_{l=0}^{n}  C_n^l
\frac{D_l^{(m)}({\bf d}^m)}{2^l} (s-\xi ({\bf d}^m))^{n-l}
$$
where
$D_l^{(m)}({\bf d}^m) \equiv 2^l B_l^{(m)}(\xi ({\bf d}^m)| {\bf d}^m)$ and
$D_{2k+1}^{(m)}({\bf d}^m) = 0$ and obtain
\begin{equation}
V_1(s,{\bf d}^m) =
\frac{1}{(m-1)!\pi({\bf d}^m)}
\sum_{l=0}^{m-1}  C_{m-1}^l s^{m-1-l}
\frac{D_{l}^{(m)}({\bf d}^m)}{2^l}.
\label{V_1f}
\end{equation}
The quantities  $D_n^{(m)}({\bf d}^m)$ can be calculated using a recursive
relation
\begin{equation}
D_n^{(m)}({\bf d}^m)=
\sum_{l=0}^{n}  C_{n}^l d_m^{l} D_l
D_{n-l}^{(m-1)}({\bf d}^{m-1}),
\label{recursD_n^m}
\end{equation}
where $D_l$ is expressed through the value of Bernoulli polynomial of
order $l$ at fixed value of argument.
$$
D_l = 2^l B_l(1/2).
$$
Instead of (\ref{recursD_n^m}) one may use more symmetric form
\begin{equation}
D_n^{(m)}({\bf d}^m)=
\sum^n C_{n}^{{\bf r}} \prod_{i=1}^{m} d_{i}^{r_i} D_{r_i},
\label{symmD_n^m}
\end{equation}
where
$$
C_{n}^{{\bf r}} = \frac{n!}{\prod_{i=1}^{m} r_i!}
$$
is a multinomial coefficient and
summation in  (\ref{symmD_n^m}) is performed over all $r_i$ such that
$$\sum_{i=1}^m r_i = n.$$
The resulting symmetric expression for the polynomial part of the shifted
restriction function is
\begin{equation}
V_1(s,{\bf d}^m) =
\frac{1}{(m-1)!\pi({\bf d}^m)}
\sum_{l=0}^{m-1}  C_{m-1}^l s^{m-1-l}
\sum^l C_{l}^{{\bf r}} \prod_{i=1}^{m} d_{i}^{r_i} B_{r_i}(1/2).
\label{V_1symm}
\end{equation}
The inner sum in the above expression can be rewritten using the symbolic
notation accepted in theory of Bernoulli polynomials (see, for example,
\cite{Norlund}):
\begin{equation}
\sum^l C_{l}^{{\bf r}} \prod_{i=1}^{m} d_{i}^{r_i} B_{r_i}(1/2)
\equiv
\left(\sum_{i=1}^m d_i \; {}^i\! B(1/2)\right)^{l},
\label{sumsymb}
\end{equation}
where powers $r_i$ of ${}^i\! B(1/2)$ are converted into indices
$$
{}^i \! B^{r_i}(1/2) \Rightarrow B_{r_i}(1/2).
$$
Using this notation one immediately arrives at the compact expression of
$V_1(s,{\bf d}^m)$
\begin{equation}
V_1(s,{\bf d}^m) =
\frac{1}{(m-1)!\pi({\bf d}^m)}
\left(s+\sum_{i=1}^m d_i \; {}^i\! B(1/2)\right)^{m-1}.
\label{V_1symmsymb}
\end{equation}
Finally, using (\ref{transform1}), we return to the expression for
$W_1(s,{\bf d}^m)$
\begin{equation}
W_1(s,{\bf d}^m) =
\frac{1}{(m-1)!\pi({\bf d}^m)}
\left(s+\sum_{i=1}^m d_i \left[1/2 + \; {}^i\! B(1/2)\right]\right)^{m-1}.
\label{W_1symmsymb}
\end{equation}
It is easily checked that the expression (\ref{V_1symmsymb}) verifies the
general recursive relation for the shifted restriction function (which is valid also
for its polynomial part):
\begin{equation}
V(s,{\bf d}^m) - V(s-d_m,{\bf d}^m) = V(s-\frac{d_m}{2},{\bf d}^{m-1})\;.
\label{mainrec}
\end{equation}


Introducing the power expansion of $V_1$ in the form
$$
V_1(s,{\bf d}^m) = \sum_{j=1}^{m} R^m_j s^{m-j},
$$
one may write for the coefficients $R^m_j$
\begin{equation}
R^m_j =
\frac{C_{m-1}^{j-1}}{(m-1)!\pi({\bf d}^m)}
\left(\sum_{i=1}^m d_i \; {}^i\! B(1/2)\right)^{j-1}.
\label{R_j^m}
\end{equation}
The recursion relation (\ref{recursD_n^m})
is equivalent to the following recursion ($1 \leq j<m$):
\begin{equation}
R^m_j =
\frac{1}{m-j} \sum_{l=0}^{j-1} d_m^{l-1} C^l_{m-1-j+l} B_l(1/2) R^{m-1}_{j-l}.
\label{recursR_j^m}
\end{equation}
It should be noted an important distinction between two last expressions
-- the formula (\ref{R_j^m}) is the explicit expression for the polynomial part
of the shifted restriction function, while the recursion presents an imcomplete
procedure (note that $j$ in  (\ref{recursR_j^m}) cannot be set equal to $m$,
so that $R^m_m$ is not defined in the framework of this procedure).

The recursion
(\ref{recursR_j^m}) is a consequence of (\ref{mainrec}),
which repeated usage leads to more general recursive
relation
\begin{equation}
V(s+\tau_m,{\bf d}^{m})=V(s,{\bf d}^{m}) +
\sum_{p=0}^{\delta_m-1}
V(s+\tau_m-\lambda_p\cdot d_m,{\bf d}^{m-1}),
\label{mainrecrepeat}
\end{equation}
where
$$
\lambda_{p} = p + 1/2,\ \ \ \ \ \delta_m = \tau_m/d_m,
$$
and
$$
\tau_m \equiv \tau({\bf d}^m) = \mbox{LCM}({\bf d}^m),
$$
where $\mbox{LCM}({\bf d}^m)$ denotes a least common multiple of the set
${\bf d}^m$.
The function $V(s,{\bf d}^{m})$ might also be written in a "polynomial" form
\begin{equation}
V(s,{\bf d}^{m})=\sum_{j=1}^{m}R^{m}_j(s) s^{m-j},
\label{31}
\end{equation}
which leads to a
recursive formula for  $\tau_m$-periodic function $R^{m}_{j}(s)$ for
$1 \le j < m$ (see \cite{FelRub})
\begin{equation}
R^{m}_j(s) =\frac{1}{m-j}\cdot
\sum_{l=0}^{j-1} d_m^{l-1} C^l_{m-1-j+l} \cdot \delta_m^{l-1}
\sum_{p=0}^{\delta_m-1}
B_l (1-\frac{\lambda_p d_m}{\tau_m})\cdot
R^{m-1}_{j-l}(s-\lambda_p\cdot d_m)\;.
\label{knf}
\end{equation}
Introducing the shift operator
\begin{equation}
{\bf S}(s,\Delta): \{{\bf S}(s,\Delta) f(s) = f(s-\Delta) \},
\label{shiftoper}
\end{equation}
we rewrite the above recursion in the form
\begin{equation}
R^{m}_j(s) =\frac{1}{m-j}\cdot
\sum_{l=0}^{j-1} d_m^{l-1} C^l_{m-1-j+l} \cdot
\left[\delta_m^{l-1}
\sum_{p=0}^{\delta_m-1}
B_l (1-\frac{\lambda_p d_m}{\tau_m}) {\bf S}(s,\lambda_p d_m)
\right] \cdot
R^{m-1}_{j-l}(s)\;.
\label{knfshift}
\end{equation}
Comparison of (\ref{knfshift}) with (\ref{recursR_j^m}) suggests
that a replacement
\begin{equation}
{}^i \! B_{r_i} (1/2) \Longrightarrow  {}^i \! {\bf B}_{r_i} =
(\tau_i/d_i)^{r_i-1}
\sum_{p_i=0}^{\tau_i/d_i-1}
B_{r_i} (1-\frac{\lambda_{p_i} d_i}{\tau_i}) {\bf S}(s,\lambda_{p_i} d_i)
\label{Breplace}
\end{equation}
may be useful in transition from the formulas for the polynomial part of the
partition function to those of for the function itself.

Setting in (\ref{knfshift}) all $R^m_j$ independent of $s$ and using the
multiplication theorem for the Bernoulli polynomials \cite{bat53}
$$
\sum_{r=0}^{m-1} B_n(x+\frac{r}{m}) = m^{-(n-1)} B_n(mx),
$$
we immediately reproduce (\ref{recursR_j^m}).
The recursion (\ref{knf}) fails to produce
$R^{m}_m(s)$, nevertheless partial information about it can be extracted.
It is useful to separate $R^{m}_m(s)$ into two terms
\begin{equation}
R^m_{m}(s)={\cal R}^m_{m}(s)+r^m_{m}(s),
\label{Rmmseparate}
\end{equation}
where
\begin{equation}
{\cal R}^m_m(s)=
\sum_{l=1}^{m-1} \frac{\tau_m^{l-1}}{l} \sum_{p=0}^{\delta_m-1}
B_l (1-\frac{\lambda_p d_m}{\tau_m})\cdot
R^{m-1}_{m-l}(s-\lambda_p\cdot d_m)\;,
\label{knf11}
\end{equation}
is $\tau_{m-1}$-periodic function, and
the other term, $d_m$-periodic function,
$$
r^m_{m}(s) = r^m_{m}(s - d_m)
$$
remains indeterminate.
In \cite{FelRub} $r^m_m(s)$ is found by applying the method of indeterminate
coefficients using knowledge of partition function zeroes. It can be checked
using (\ref{31}) that the partition function may be written also as
\begin{equation}
V(s,{\bf d}^m) = r^m_m(s) +
\sum_{l=1}^{m-1} \frac{\tau_m^{l-1}}{l}
\sum_{p=0}^{\delta_m-1}
B_{l}\left(\frac{s+\lambda_p d_m}{\tau_m}\right)
R^{m-1}_{m-l}(s+\lambda_p d_m).
\label{Vcompact}
\end{equation}
The corresponding expression for the polynomial part of the partition function
reads
\begin{equation}
V_1(s,{\bf d}^m) = r^m_m +
\sum_{l=1}^{m-1} \frac{d_m^{l-1}}{l} \;
B_{l}\left(\frac{1}{2}+\frac{s}{d_m}\right)
R^{m-1}_{m-l}.
\label{V1compact}
\end{equation}
Setting in the above expression $s=0$ we obtain polynomial analog of
(\ref{Rmmseparate})
\begin{equation}
R^m_{m}=r^m_{m} +
\sum_{l=1}^{m-1} \frac{d_m^{l-1}}{l}
B_l (\frac{1}{2}) R^{m-1}_{m-l}.
\label{Rmm1sep}
\end{equation}
Substituting into it the general expression (\ref{R_j^m}) we determine
$r^m_{m}$
\begin{equation}
r^m_m =
\frac{1}{(m-1)!\;\pi({\bf d}^m)}
\left(\sum_{i=1}^{m-1} d_i \; {}^i\! B(1/2)\right)^{m-1}
.
\label{r_m^m}
\end{equation}


\section{Calculation of $r^m_m(s)$}

Our goal is to find yet indeterminate $r^m_m(s)$ using (\ref{r_m^m})
and recursive relation (\ref{mainrec}) for the partition function and its
polynomial part and its corrolaries (\ref{recursR_j^m}) and (\ref{knf}).

In a simplest case $m=1$ it is easily checked that the
shifted partition function $V(s,\{d_1\})$ is a $d_1$-periodic function given
by
\begin{equation}
V(s,\{d_1\}) = R^1_1(s) = r^1_1(s) = \Psi_{d_1}(s-\frac{d_1}{2})=
\sum_{p_1=0}^{\tau_1/d_1-1} B_0(\frac{\lambda_{p_1} d_1}{\tau_1})
\Psi_{d_1}(s-\lambda_{p_1} d_1),
\label{V1}
\end{equation}
Here $\tau_1 \equiv d_1$, so that $p_1$ takes only zero value; we also
use parity property of Bernoulli polynomials $B_0(1-x)=B_0(x)
\equiv 1$.  The periodic function $\Psi_{d_1}(s)$ is defined as
a sum of prime roots of
unit of degree $d_1$:
$$
\Psi_{d_1}(s) =
\frac{1}{d_1}\sum_{k=0}^{d_1-1} \exp \left(\frac{2 \pi i k s}{d_1}\right) =
\left\{ \begin{array}{ll}
         1\;\;,   & \mbox{$s=0 \pmod{d_1}$}  \\
         0\;\;,   & \mbox{$s\neq 0 \pmod{d_1}$}
                           \end{array}\right.
$$
The polynomial part of this function
$$
V_1(s,\{d_1\}) = \frac{1}{d_1},
$$
this also follows from (\ref{V_1symmsymb}) for $m=1$.

Consider $m=2$, and find $R^2_1(s)$ using the recursive relation (\ref{knf})
for the set $\{d_1,d_2\}$
\begin{eqnarray}
R^2_1(s) & = & \frac{1}{\tau_2} \sum_{p_2=0}^{\tau_2/d_2-1}
B_0(1-\frac{\lambda_{p_2} d_2}{\tau_2}) R^1_1(s - \lambda_{p_2} d_2) =
\nonumber \\
&&
\frac{1}{\tau_2} \sum_{p_2=0}^{\tau_2/d_2-1}
B_0(1-\frac{\lambda_{p_2} d_2}{\tau_2})
\sum_{p_1=0}^{\tau_1/d_1-1} B_0(1-\frac{\lambda_{p_1} d_1}{\tau_1})
\Psi_{d_1}(s-\lambda_{p_2} d_2-\lambda_{p_1} d_1).
\label{R21}
\end{eqnarray}
The symmetry of the problem w.r.t. the permutations of the set elements implies
that one can apply an interchange
$d_1 \leftrightarrow d_2$ in order to produce another valid form of $R^2_1(s)$.
The corresponding polynomial part is found as
$$
R^2_1 = \frac{{}^1 \! B_0(1/2)\; {}^2 \! B_0(1/2)}{d_1 d_2} \;.
$$
It is clear that (\ref{R21}) can be produced from the above as
follows -- we replace the fraction $1/d_1$ by its counterpart $\Psi_{d_1}(s -
d_1/2)$ and then use the replacement (\ref{Breplace}). The another possible
form is found by the replacement $1/d_2 \Rightarrow \Psi_{d_2}(s -
d_2/2)$ and application of (\ref{Breplace}) to it.
The free term $R^2_2(s)$ is constructed in two steps. At first we find
${\cal R}^2_2(s)$ assuming $d_1$ the first element of the set. Then
(\ref{knf11}) gives
\begin{eqnarray}
{\cal R}^2_2(s) & = & \sum_{p_2=0}^{\tau_2/d_2-1}
B_1(1-\frac{\lambda_{p_2} d_2}{\tau_2}) R^1_1(s - \lambda_{p_2} d_2) =
\nonumber \\
&&
\sum_{p_2=0}^{\tau_2/d_2-1}
B_1(1-\frac{\lambda_{p_2} d_2}{\tau_2})
\sum_{p_1=0}^{\tau_1/d_1-1} B_0(1-\frac{\lambda_{p_1} d_1}{\tau_1})
\Psi_{d_1}(s-\lambda_{p_2} d_2-\lambda_{p_1} d_1).
\label{R22a}
\end{eqnarray}
The other part -- $d_2$-periodic function $r^2_2(s)$ may be produced in this
case by the interchange
$d_1 \leftrightarrow d_2$ in the above expression, which corresponds to choice
of $d_2$ as the first element of the set.
\begin{eqnarray}
r^2_2(s) & = & \sum_{p_1=0}^{\tau_2/d_1-1}
B_1(1-\frac{\lambda_{p_1} d_1}{\tau_2}) R^1_1(s - \lambda_{p_2} d_2) =
\nonumber \\
&&
\sum_{p_1=0}^{\tau_2/d_1-1}
B_1(1-\frac{\lambda_{p_1} d_1}{\tau_2})
\sum_{p_2=0}^{\tau_1/d_2-1} B_0(1-\frac{\lambda_{p_2} d_2}{\tau_1})
\Psi_{d_2}(s-\lambda_{p_2} d_2-\lambda_{p_1} d_1).
\label{r22}
\end{eqnarray}
The corresponding polynomial parts reads
$$
r^2_2 = \frac{{}^1 \! B_1(1/2)\; d_1}{d_1 d_2}
$$
and is equal to zero.
Nevertheless, we may use  the above expression in order to produce
(\ref{r22}) -- we start from the replacement $1/d_2 \Rightarrow \Psi_{d_2}(s -
d_2/2)$ and then apply (\ref{Breplace}) to arrive at $d_2$-periodic function.

We conjecture the following formula for $d_m$-periodic function
$r^m_m(s)$ using its polynomial part (\ref{r_m^m}) which we rewrite in the
expanded form
\begin{equation}
r^m_m = \frac{1}{d_m} \cdot
\frac{1}{(m-1)!}
\sum^{m-1} C_{m-1}^{{\bf r}}
\prod_{i=1}^{m} d_{i}^{\;r_i-1} \; {}^i \! B_{r_i}(1/2).
\label{r_m^mexpand}
\end{equation}

We start with
$1/d_m \Rightarrow \Psi_{d_m}(s - d_m/2)$ securing the proper
periodicity. Then the replacements (\ref{Breplace}) lead to the form
\begin{equation}
r^m_m(s) =
\frac{1}{(m-1)!}
\sum^{m-1} C_{m-1}^{{\bf r}}
\prod_{i=1}^{m-1} d_{i}^{\;r_i-1} \; {}^i {\bf B}_{r_i}
\Psi_{d_m}(s - d_m/2).
\label{r_m^m(s)}
\end{equation}
Using the explicit form of the operator ${}^i {\bf B}_{r_i}$ we arrive at
\begin{equation}
r^m_m(s) =
\frac{1}{(m-1)!}
\sum^{m-1}_{{\bf r}} C_{m-1}^{{\bf r}}
\prod_{i=1}^{m-1} \tau_{i}^{r_i-1}
\sum_{p_i=0}^{\tau_i/d_i-1}
 B_{r_i}(1-\frac{\lambda_{p_i} d_i}{\tau_i})
\Psi_{d_m}(s - \sum_{i=1}^{m} \lambda_{p_i} d_i).
\label{r_m^m(s)sums}
\end{equation}
It should be underlined here that the
values of "periods" $\tau_i$ depends on the order of the elements
in the set ${\bf d}^m$,
and in the above formula we have $\tau_1=\mbox{LCM}(d_m,d_1),
\tau_2=\mbox{LCM}(d_m,d_1,d_2),\ldots,
\tau_{m-1}=\mbox{LCM}({\bf d}^m), \tau_m = d_m$, so that $p_m \equiv 0$.
The last formula can be rewritten in the symbolic form
\begin{equation}
r^m_m(s) =
\frac{1}{(m-1)!}
\prod_{i=1}^{m-1} \tau_{i}^{-1}
\left[
\sum^{m-1}_{i=1} \tau_{i}
\left\{
\sum_{p_i=0}^{\tau_i/d_i-1}
{}^i \! B(1-\frac{\lambda_{p_i} d_i}{\tau_i})
\right \}
\right]^{m-1}
\Psi_{d_m}(s - \sum_{i=1}^{m} \lambda_{p_i} d_i).
\label{r_m^m(s)symb}
\end{equation}
Hence, combination of the last expression with the recursion (\ref{Vcompact})
provides a new procedure for calculation of the restricted partition
function  $V(s,{\bf d}^m)$.


\section{Explicit formula for restricted partition function}

The same approach can be useful in order to produce a formula for
$V(s,{\bf d}^m)$, starting from its polynomial part
(\ref{V_1symm}) which can be written as
\begin{equation}
V_1(s,{\bf d}^m) =
\frac{1}{(m-1)!}
\left(\prod_{i=1}^{m} d_{i}^{-1} \right)
\sum_{l=0}^{m-1}  C_{m-1}^l s^{m-1-l}
\left(\sum_{i=1}^m d_i \; {}^i\! B(1/2)\right)^{l}.
\label{V_1symmnew}
\end{equation}
It is clear that there exist several equivalent forms of partition function
(like two possible forms of $R^2_1(s)$ discussed above), and
symmetry considerations can help in selection of the most symmetric form. In the
discussed case the form is chosen by a selection of the factor $1/d_i$ for the
replacement $1/d_i \Rightarrow \Psi_{d_i}(s-d_i/2)$
and application of (\ref{Breplace}) to the result.

In general situation we have $m$ choices of $1/d_i$ to start with, and because
the result doesn't depend of such choice, it is natural to seek
an expression symmetrized w.r.t. the starting element $d_i$. The latter problem
is eqivalent to the problem of presentation of the symmetric polynomial
$$
\left(\sum_{i=1}^m d_i \right)^{l}=
\sum^l C_{l}^{{\bf r}} \prod_{i=1}^{m} d_{i}^{r_i}
$$
as a sum of $m$ {\it symmetric} polynomials, each of them missing only one of
terms $d_i$. Introducing a function $z({\bf r})$ counting number of zero
components of the vector ${\bf r}$, one can write the above expression as
\begin{equation}
\sum_{i=1}^m
\sum^l_{\bf r} \frac{1}{z({\bf r})} C_{l}^{{\bf r}}
\prod_{n=1 \atop n \ne i}^{m} d_{n}^{\;r_n}.
\label{symmpol1}
\end{equation}
Using this result we have
\begin{equation}
\left( \prod_{i=1}^{m} d_{i}^{-1} \right)
\left(\sum_{i=1}^m d_i \; {}^i\! B(1/2)\right)^{l} =
\sum_{i=1}^m  \frac{1}{d_i}
\sum^l_{\bf r} \frac{1}{z({\bf r})} C_{l}^{{\bf r}}
\prod_{n=1 \atop n \ne i}^{m} d_{n}^{\;r_n-1} B_{r_n}(1/2).
\label{symmpol2}
\end{equation}
The polynomial part (\ref{V_1symmnew})
suitable for conversion into $V(s,{\bf d}^m)$ has the form
\begin{equation}
V_1(s,{\bf d}^m) =
\frac{1}{(m-1)!}
\sum_{l=0}^{m-1}  C_{m-1}^l s^{m-1-l}
\sum_{i=1}^m  \frac{1}{d_i}
\sum^l_{\bf r} \frac{1}{z({\bf r})} C_{l}^{{\bf r}}
\prod_{n=1 \atop n \ne i}^{m} d_{n}^{\;r_n-1} B_{r_n}(1/2).
\label{V_1symmnewer}
\end{equation}
The explicit expression for the restricted
partition function reads
\begin{equation}
V(s,{\bf d}^m) =
\frac{1}{(m-1)!}
\sum_{l=0}^{m-1}  C_{m-1}^l s^{m-1-l}
\sum_{i=1}^m
\left[
\sum^l_{\bf r} \frac{1}{z({\bf r})} C_{l}^{{\bf r}}
\prod_{n=1 \atop n \ne i}^{m} d_{n}^{\;r_n-1} \; {}^{n} {\bf B}_{r_n}
\right]
\Psi_{d_i}(s-d_i/2).
\label{Vsymm}
\end{equation}
Using the actual form of the operators
${}^{n} {\bf B}_{r_n}$ (\ref{Breplace}) it is rewritten as
\begin{eqnarray}
V(s,{\bf d}^m) & = &
\frac{1}{(m-1)!}
\sum_{l=0}^{m-1}  C_{m-1}^l s^{m-1-l} \\
&&
\sum_{i=1}^m
\left[
\sum^l_{\bf r} \frac{1}{z({\bf r})} C_{l}^{{\bf r}}
\prod_{n=1 \atop n \ne i}^{m} \tau_{n,i}^{\;r_n-1}
\sum_{p_n=0}^{\tau_{n,i}/d_n-1} B_{r_n}(1-\frac{\lambda_{p_n} d_n}{\tau_{n,i}})
\right]
\Psi_{d_i}(s-\sum_{n=1}^m \lambda_{p_n} d_n).   \nonumber
\label{Vsymmfinal}
\end{eqnarray}
It should be noted here that the values of "periods" $\tau$ depends on selected
value $d_i$ what is reflected by an additional subscript $i$:
\\
\\
\begin{tabular}{ccccc}
$\tau_{1,1} = d_1 $&
$\tau_{2,1}=\mbox{LCM}(d_1,d_2) $&
$\tau_{3,1}=\mbox{LCM}(d_1,d_2,d_3) $&
$\ldots $&
$\tau_{m,1}=\mbox{LCM}({\bf d}^m) $\\
$\tau_{1,2}=\mbox{LCM}(d_1,d_2) $&
$\tau_{2,2}=d_2$ &
$\tau_{3,2}=\mbox{LCM}(d_1,d_2,d_3)$ &
$\ldots$ &
$\tau_{m,2}=\mbox{LCM}({\bf d}^m)$ \\
$\tau_{1,3}=\mbox{LCM}(d_1,d_3) $ &
$\tau_{2,3}=\mbox{LCM}(d_1,d_2,d_3)  $ &
$\tau_{3,3}=d_3 $ &
$\ldots $ &
$ \tau_{m,3}=\mbox{LCM}({\bf d}^m) $\\
$\ldots $ &
$\ldots $ &
$\ldots $ &
$\ldots $ &
$\ldots $ \\
$\tau_{1,m}=\mbox{LCM}(d_1,d_m) $ &
$\tau_{2,m}=\mbox{LCM}(d_1,d_2,d_m) $ &
$\tau_{3,m}=\mbox{LCM}(d_1,d_2,d_3,d_m) $  &
$\ldots $ &
$\tau_{m,m}=d_m $
\end{tabular}




\begin{thebibliography}{99}
\bibitem{Euler}
L. Euler, {\it Introductio in Analysin Infinitorum},
Vol. I, Lausanne (1748).
\bibitem{GAndrews}
G. E. Andrews, {\it The Theory of Partitions},
\\Encyclopedia of Mathematics and
its Applications, V.2, Addison-Wesley (1976).
\bibitem{Cayley} A. Cayley, {\it Researches on the partitions of numbers},
\\
Phyl. Trans. Royal Soc. {\bf 145} p.127-140 (1855),
\\reprinted in
{\it Coll. Papers}, Vol. 2, p.235-249, Cambridge University Press, Cambridge
(1889).
\bibitem{Sylv0} J.J. Sylvester, {\it A Constructive Theory of Partitions,
Arranged in Three acts, an Interact and an Exodus}, \\
American Journal of Mathematics {\bf 5}, p.251-300 (1882).
\bibitem{Sylv1} J.J. Sylvester, {\it On the Partition
of Numbers}, \\
 Quarterly Journal of Mathematics {\bf I}, p.141-152 (1857).
\bibitem{Sylv2} J.J. Sylvester, {\it On Subinvariants, i.e. Semi-invariants
to Binary Quantics of an Unlimited Order. With an Excursus on Rational
Fractions and Partitions},\\
 American Journal of Mathematics {\bf 5}, p.79-136 (1882).
\bibitem{Herschel} J.F.W. Herschel, {\it On Circulating Functions, and on the
Integration of a Class of Equations od Finite Differences into which They Enter
as Coefficients},\\
 Phyl. Trans. Royal Soc. {\bf 108}, p.144-168 (1818).
\bibitem{Gupta}  H. Gupta, {\it Tables of Partitions}, \\
Royal Society
Mathematical Tables, v.4, Cambridge University Press (1958).
\bibitem{FelRub} L.G. Fel and B.Y. Rubinstein,
{\it Sylvester Waves in the Coxeter Groups},
\\submitted to Ramanujan Journal (2001);
\\math-nt/0005174, deposited in the server {\tt http://xxx.arXiv.org} (2000).
\bibitem{bat53}    H.Bateman and A.Erdel\'yi, {\it Higher Transcendental
Functions}, {\bf 1},
\\ McGraw-Hill Book Co,Inc., N.Y., (1953).
\bibitem{Norlund}  N.E. N\"orlund {\it Vorlesungen \"uber Differenzenrechnung},
\\Verlag von Julius Springer, Berlin (1924).
\end{thebibliography}
\end{document}